\documentclass[11pt]{article}

\usepackage{mathpazo}   
\usepackage{tikz}

\usepackage{nameref}
\usepackage{empheq}
\usepackage{comment}
\usepackage[shortlabels,inline]{enumitem}
\setlist[enumerate]{nosep}
\usepackage[colorlinks=true,
linkcolor=refkey,
urlcolor=lblue,
citecolor=red]{hyperref}
\usepackage[doc,wmm,hhb]{optional}
\usepackage{xcolor}

\usepackage{float}
\usepackage{soul}
\usepackage{graphicx}
\definecolor{labelkey}{rgb}{0,0.08,0.45}
\definecolor{refkey}{rgb}{0,0.6,0.0}
\definecolor{Brown}{rgb}{0.45,0.0,0.05}
\definecolor{lime}{rgb}{0.00,0.8,0.0}
\definecolor{lblue}{rgb}{0.5,0.5,0.99}
\definecolor{OliveGreen}{rgb}{0,0.6,0}
\definecolor{tyrianpurple}{rgb}{0.4, 0.01, 0.24}

\colorlet{hlcyan}{cyan!30}

\usepackage{stmaryrd}
\usepackage{amssymb}

\makeatletter
\def\namedlabel#1#2{\begingroup
   \def\@currentlabel{#2}%
   \label{#1}\endgroup
}
\makeatother

\newcommand{\seppthree}{\setlength{\itemsep}{-3pt}}

\usepackage[margin=0.92in,footskip=0.25in]{geometry}
\providecommand{\siff}{\Leftrightarrow}

\newcommand{\thalb}{\ensuremath{\tfrac{1}{2}}}
\newcommand{\menge}[2]{\big\{{#1}~\big |~{#2}\big\}}

\newcommand{\fenv}[1]%
{\ensuremath{\,\overrightarrow{\operatorname{env}}_{#1}}}
\newcommand{\benv}[1]%
{\ensuremath{\,\overleftarrow{\operatorname{env}}_{#1}}}

\newcommand{\scal}[2]{\left\langle{#1},{#2}  \right\rangle}

\newcommand{\tscal}[2]{\langle{#1},{#2}  \rangle}

\newcommand{\zeroun}{\ensuremath{\left]0,1\right[}}
\newcommand{\RR}{\ensuremath{\mathbb R}}
\newcommand{\RP}{\ensuremath{\mathbb{R}_+}}
\newcommand{\RPP}{\ensuremath{\mathbb{R}_{++}}}

\newcommand{\ran}{\ensuremath{{\operatorname{ran}}\,}}

\newcommand{\Fix}{\ensuremath{\operatorname{Fix}}}
\newcommand{\Id}{\ensuremath{\operatorname{Id}}}

\newcommand{\pinf}{\ensuremath{+\infty}}

{\begin{list}{}{%
\settowidth{\labelwidth}{\textrm{#1~}}%
\setlength{\leftmargin}{\labelwidth+\labelsep}}}
{\end{list}}
\usepackage{amsthm}
\makeatletter
\def\th@plain{%
	\thm@notefont{}
	\itshape 
}
\def\th@definition{%
	\thm@notefont{}
	\normalfont 
}
\makeatother
\usepackage[capitalize,nameinlink]{cleveref}
\crefname{equation}{}{equations}
\crefname{chapter}{Appendix}{chapters}
\crefname{item}{}{items}
\crefname{enumi}{}{}
\newtheorem{theorem}{Theorem}[section]
\newtheorem{lemma}[theorem]{Lemma}

\newtheorem{corollary}[theorem]{Corollary}

\newtheorem{proposition}[theorem]{Proposition}

\newtheorem{definition}[theorem]{Definition}

\newtheorem{example}[theorem]{Example}
\newtheorem{ex}[theorem]{Example}
\newtheorem{fact}[theorem]{Fact}
\newtheorem{remark}[theorem]{Remark}

\providecommand{\RR}{\mathbb{R}}

\providecommand{\ran}{\operatorname{ran}}

\providecommand{\Id}{\operatorname{{ Id}}}

\providecommand{\ran}{\operatorname{ran}}

\providecommand{\Id}{\operatorname{Id}}

\providecommand{\spn}{\operatorname{span}}

\newcommand{\cran}{\ensuremath{\overline{\operatorname{ran}}\,}}

\providecommand{\RR}{\mathbb{R}}

\definecolor{myblue}{rgb}{0.9,0.9,0.98}

  \newcommand*\mybluebox[1]{%
    \colorbox{myblue}{\hspace{1em}#1\hspace{1em}}}

\allowdisplaybreaks 

\newcommand*{\tran}{^{\mkern-1.5mu\mathsf{T}}}

\begin{document}

\setlength{\abovedisplayskip}{8pt}
\setlength{\belowdisplayskip}{8pt}	
	

%

\author{
Heinz H.\ Bauschke\thanks{
Mathematics, University
of British Columbia,
Kelowna, B.C.\ V1V~1V7, Canada. E-mail:
\texttt{heinz.bauschke@ubc.ca}.},~~~
Theo Bendit\thanks{
Mathematics, University
of British Columbia,
Kelowna, B.C.\ V1V~1V7, Canada. E-mail:
\texttt{theo.bendit@gmail.com}.},
~~~and~
Walaa M.\ Moursi\thanks{
Department of Combinatorics and Optimization, 
University of Waterloo,
Waterloo, Ontario N2L~3G1, Canada.
  E-mail: \texttt{walaa.moursi@uwaterloo.ca}.}
}

\title{\textsf{
How averaged is the composition of two linear projections?
}
}

\date{March 24, 2023}

\maketitle

\begin{abstract}
Projection operators are fundamental algorithmic operators in Analysis and Optimization. It is well known that these operators are firmly nonexpansive; 
however, their composition is generally only averaged and no longer firmly nonexpansive. In this note, we introduce the modulus of averagedness and provide 
an exact result for the composition of two linear projection operators. 
As a consequence, we deduce that 
the Ogura--Yamada bound for the modulus of the composition is sharp. 
\end{abstract}
{ 

}

\section{Introduction}

Throughout, we assume that 
\begin{empheq}[box=\mybluebox]{equation}
\text{$X$ is
a real Hilbert space with inner product 
$\scal{\cdot}{\cdot}\colon X\times X\to\RR$, }
\end{empheq}
and induced norm $\|\cdot\|$.

Let $T\colon X\to X$ be nonexpansive. 
Recall that $T$ is \emph{$\kappa$-averaged}\footnote{Usually, one excludes the cases $\kappa=0$ and $\kappa=1$ in the study of averaged operators, but it is very convenient in this paper to allow for this case.}, if 
$T$ can be represented as 
$T=(1-\kappa)\Id+\kappa N$, where $N$ is nonexpansive and 
$\kappa\in[0,1]$.
Note that the identity operator, written $\Id$, is the only $0$-averaged operator. 
It is clear that if $T$ averaged, then it is nonexpansive. 
Various characterizations of $\kappa$-averagedness are available including 
(see \cite[Proposition~4.35]{BC2017}) 
\begin{equation}
\label{e:230228a}
(\forall x\in X)(\forall y\in X)\quad
\|Tx-Ty\|^2\leq\|x-y\|^2 - \frac{1-\kappa}{\kappa}\|(\Id-T)x-(\Id-T)y\|^2, 
\end{equation}
which does require $\kappa>0$, and also 
\begin{equation}
\label{e:230228aa}
(\forall x\in X)(\forall y\in X)\quad
\|Tx-Ty\|^2 + (1-2\kappa)\|x-y\|^2\leq 2(1-\kappa)\scal{x-y}{Tx-Ty}.
\end{equation}
In the analysis of algorithms, one is often interested in making 
the estimate \cref{e:230228a} as sharp as possible. 
Note the the difference in the right side of \cref{e:230228a} is smallest if the subtrahend is as large as possible, 
i.e., if $\kappa$ is as small as possible. 
The operator $T$ is called \emph{averaged} if it is $\kappa$-averaged for some
$\kappa\in\left[0,1\right[$. Clearly, this is the interesting case because only 
then does 
\cref{e:230228a} provide a nontrivial estimate. 
In light of this the following definition is now quite natural:

\begin{definition}[modulus of averagedness]
\label{d:moa}
Let $T\colon X\to X$ be nonexpansive\footnote{We assume for convenience throughout the paper that the operators have full domain which is the case in all algorithmic applications we are aware of. One could obviously generalize this notion to allow for operators whose domains are proper subsets of $X$.}.
Then we define the \emph{modulus of averagedness} by 
\begin{equation}
\label{e:moa}
\kappa(T) := \min\menge{\kappa\in[0,1]}{\text{$T$ is $\kappa$-averaged}}. 
\end{equation}
\end{definition}

\begin{remark}
Note that in \cref{e:moa} we deliberately wrote ``$\min$'' and not ``$\inf$''.
Indeed, the characterization of $\kappa$-averaged in \cref{e:230228aa} makes this clear as both the left and right hand side are continuous functions in $\kappa$. 
\end{remark}

Prime examples of averaged operators are projections and more generally proximal mappings or even resolvents. These operators are \emph{firmly nonexpansive}, i.e., $\thalb$-averaged. 
However, the class of firmly nonexpansive operators is
\emph{not} closed under compositions 
(see, e.g., \cref{ex:wow} below).
This defect is remedied by 
considering averaged operators which are closed under convex combinations and 
compositions (see, e.g., \cite{BC2017}). 
More quantitatively, and formulated here using the modulus of averagedness, 
Ogura and Yamada obtained the following:

\begin{fact}[Ogura--Yamada] {(See \cite[Theorem~3]{OY}.)}
\label{f:OY}
Let $T_1\colon X \to X$ and 
$T_2\colon X \to X$ be averaged, 
and let $\lambda\in[0,1]$. 
Then 
\begin{equation}
\label{e:OY}
\kappa(T_2T_1)\leq \frac{\kappa(T_1)+\kappa(T_2)-2\kappa(T_1)\kappa(T_2)}{1-\kappa(T_1)\kappa(T_2)}
\end{equation}
and
$\kappa\big((1-\lambda)T_1+\lambda T_2\big)\leq (1-\lambda)\kappa(T_1)+\lambda\kappa(T_2)$. 
\end{fact}

We are now ready to explain the goal of this note.
We will study the Ogura--Yamada bound \cref{e:OY} 
in the context of the composition of two 
linear projections (one of which may even be relaxed) onto two
closed linear subspaces $U,V$ of $X$. Our main result (see \cref{t:main} below) implies that 
the Ogura--Yamada bound can be improved in this setting when 
the \emph{Friedrichs angle} (see \cref{e:c_F}) 
between $U$ and $V$ is positive 
(which is always the case when $X$ is finite-dimensional) and 
also that the Ogura--Yamada bound is \emph{sharp} if the Friedrichs angle is $0$.

The remainder of this paper is organized as follows.
Some auxiliary and nonlinear results are collected in \cref{sec:aux}.
The main result is developed in \cref{sec:main}.

The notation employed in this paper is fairly standard and as in, e.g.,  \cite{BC2017}.

\section{Some linear and nonlinear observations}

\label{sec:aux}

\subsection{The modulus of averagedness in the linear case}

\begin{lemma}
\label{l:moachar}
Let $T\colon X\to X$ be linear and nonexpansive, 
and suppose that $T\neq\Id$. Then
\begin{equation}
\label{e:maochar}
\kappa(T) = \sup_{z\neq Tz } \frac{\|z-Tz\|^2}{2\scal{z}{z-Tz}}, 
\end{equation}
and the quotient appearing in the supremum is well defined. 
\end{lemma}
\begin{proof}
We start by noting that the quotient appearing in \cref{e:maochar} is well defined.
Indeed, 
Cauchy--Schwarz yields
$\scal{z}{Tz}\leq\|z\|\|Tz\|\leq \|z\|^2$ because $T$ is nonexpansive. 
Thus if 
$\scal{z}{Tz}=\|z\|^2$, then the equality condition for Cauchy--Schwarz yields $z=Tz$. Putting this contrapositively, if $z\neq Tz$, then 
$\scal{z}{Tz}<\|z\|^2$ and the quotient is thus well defined. 

Now let $\kappa\in\left]0,1\right]$. 
Using the linearity of $T$, we see that 
\cref{e:230228a} turns into 
\begin{equation}
\label{e:230228b}
\|Tz\|^2\leq\|z\|^2 - \frac{1-\kappa}{\kappa}\|(\Id-T)z\|^2,
\end{equation}
where $z\in X$.
If $z=Tz$, then \cref{e:230228b} is obviously true (even for any $\kappa$).
Thus assume that $z\neq Tz$. 
Then \cref{e:230228b} is equivalent to 
\begin{equation}
\label{e:230228bb}
\frac{\|z\|^2-\|Tz\|^2}{\|z-Tz\|^2} \geq \frac{1-\kappa}{\kappa}.
\end{equation}
In turn, \cref{e:230228bb} is equivalent to 
\begin{subequations}
\begin{align}
\frac{1}{\kappa}
&\leq 1+ \frac{\|z\|^2-\|Tz\|^2}{\|z-Tz\|^2}
= \frac{\|z-Tz\|^2+\|z\|^2-\|Tz\|^2}{\|z-Tz\|^2}\\
&= 
\frac{\|z\|^2+\|Tz\|^2-2\scal{z}{Tz}+\|z\|^2-\|Tz\|^2}{\|z-Tz\|^2}\\
&= \frac{2\big(\|z\|^2-\scal{z}{Tz} \big)}{\|z-Tz\|^2}
= \frac{2\scal{z}{z-Tz}}{\|z-Tz\|^2}.
\end{align}
\end{subequations}
The proof is complete.
\end{proof}

\subsection{Underrelaxed nonexpansive linear operators}

The next result shows that the computation of the modulus of averagedness
can sometimes be reduced to a closed linear subspace of $X$.

\begin{theorem}
\label{t:230228}
Let $R\colon X\to X$ be linear and nonexpansive.
Set  $V := \Fix R$ and assume that $V\subsetneqq X$.
Let $\beta\in[0,1]$, and set 
$T := (1-\beta)\Id + \beta R$.
Then 
\begin{equation}
\kappa(T) = \beta\kappa\big(P_{V^\perp}R|_{V^\perp}\big).
\end{equation}
\end{theorem}
\begin{proof}
The conclusion is clear if $\beta=0$. 
So assume that $\beta>0$. 
Note that $\Fix T = \Fix R = \ker(\Id-R)=V$.
By \cite[Example~20.29]{BC2017}, $\Id-R$ is maximally monotone.
By \cite[Proposition~20.17]{BC2017}, $\ker(\Id-R)=\ker(\Id-R^*)$
and $\cran(\Id-R)=\cran(\Id-R^*)$.
Hence $\cran(\Id-R)=(\ker(\Id-R^*))^\perp = (\ker(\Id-R))^\perp = V^\perp$.
Thus 
\begin{equation}
V = \Fix R = \ker(\Id-R)
\;\;\text{and}\;\; \cran(\Id-R)=V^\perp.
\end{equation}
Now let $z\in X\smallsetminus \Fix T = X\smallsetminus \Fix R = X\smallsetminus V$, and decompose $z=P_Vz+P_{V^\perp}z$.
Then 
$z-Rz
=(P_Vz+P_{V^\perp}z)-(R(P_Vz+P_{V^\perp}z))
=P_Vz+P_{V^\perp}z - R(P_Vz)- RP_{V^\perp}z
=P_Vz+P_{V^\perp}z - P_Vz- RP_{V^\perp}z
= (\Id-R)P_{V^\perp}z
\in \cran(\Id-R)
= V^\perp$; to sum up,
\begin{equation}
z-Rz = (\Id-R)P_{V^\perp}z\in V^\perp 
\;\;\text{and so}\;\;
z-Rz = P_{V^\perp}z-P_{V^\perp}RP_{V^\perp}z.
\end{equation}
Recalling the definition of $T$, we obtain 
\begin{subequations}
\begin{align}
\frac{\|z-Tz\|^2}{2\scal{z}{z-Tz}}
&= \frac{\|z-(1-\beta)z-\beta Rz\|^2}{2\scal{z}{z-(1-\beta)z-\beta Rz}}
= \frac{\|\beta(z-Rz)\|^2}{2\scal{z}{\beta(z-Rz)}}\\
&= \beta\frac{\|z-Rz\|^2}{2\scal{z}{z-Rz}}
= \beta\frac{\|P_{V^\perp}z-P_{V^\perp}RP_{V^\perp}z\|^2}{2\scal{P_Vz+P_{V^\perp}z}{P_{V^\perp}z-P_{V^\perp}RP_{V^\perp}z}}\\
&= \beta\frac{\|(P_{V^\perp}z)-(P_{V^\perp}R)(P_{V^\perp}z)\|^2}{2\scal{(P_{V^\perp}z)}{(P_{V^\perp}z)-(P_{V^\perp}R)(P_{V^\perp}z)}}, 
\end{align}
\end{subequations}
which features the quotient arising in the characterization of the modulus of
averagedness of 
the operator 
$P_{V^\perp}R|_{V^\perp}\colon V^\perp\to V^\perp$ (see \cref{l:moachar}). 
\end{proof}

\begin{example}
\label{ex:sarah}
Let $V$ be a closed linear subspace of $X$, 
define the corresponding reflector by
$R_V := P_V-P_{V^\perp}$, and let $\beta\in[0,1]$. 
Then $\Fix R_V = V$.
and 
\begin{equation}
\kappa\big((1-\beta)\Id+\beta R_V \big) = 
\begin{cases}
0, &\text{if $V=X$;}\\
\beta, &\text{if $V\neq X$.}
\end{cases}
\end{equation}
\end{example}
\begin{proof}
The conclusion is clear when $\beta=0$ so we assume that $\beta>0$. 
If $V=X$, then $R_V=\Id$ and thus $(1-\beta)\Id+\beta R_V=\Id$ and the conclusion follows. Now assume that $V\neq X$ which yields $\{0\}\subsetneqq V^\perp$. 
Clearly, 
$R_V|_{V^\perp}=P_V|_{V^\perp}-P_{V^\perp}|_{V^\perp}
= -\Id_{V^\perp}$
which implies
$P_{V^\perp}R_V|_{V^\perp}=-\Id_{V^\perp}$.
Next, let $w\in V^\perp$ such that $w\notin \Fix(-\Id_{V^\perp})$, i.e., 
$w\in V^\perp\smallsetminus\{0\}$. Then 
\begin{align}
\frac{\|w-(-\Id)w\|^2}{2\scal{w}{w-(-\Id)w}}
&= 
\frac{\|2w\|^2}{2\scal{w}{2w}} = 1
\end{align}
and we deduce that 
$\kappa(P_{V^\perp}R_V|_{V^\perp})=1$. 
The conclusion now follows from \cref{t:230228}. 
\end{proof}

\subsection{Results that involve adjoint operators}

The modulus of averagedness is invariant under taking the adjoint: 

\begin{proposition}
\label{p:mod*}
Let $T\colon X\to X$ be linear and nonexpansive. 
Then $\kappa(T)=\kappa(T^*)$. 
\end{proposition}
\begin{proof}
Indeed, let $\kappa\in[0,1]$.
If $T$ is $\kappa$-averaged, then there exists
$N\colon X\to X$ nonexpansive such that 
$T=(1-\kappa)\Id+\kappa N$.
Taking adjoints yields
$T^*=(1-\kappa)\Id+\kappa N^*$.
Because $N^*$ is also nonexpansive 
(recall that $\|N\|=\|N^*\|$ by, e.g., \cite[Theorem~8.25]{Deutsch}), we see that $T^*$ is also $\kappa$-averaged.
This implies $\kappa(T^*)\leq\kappa(T)$. 
Applying this to $T^*$ instead of $T$ yields
$\kappa(T)=\kappa(T^{**})\leq\kappa(T^*)$ and we are done.
\end{proof}

\begin{lemma}
\label{l:nonexpchar}
Let $T\colon X\to X$ be linear and continuous.
Then $T$ is nonexpansive if and only if 
$\Id-T^*T$ is positive semidefinite.
\end{lemma}
\begin{proof}
Using the linearity of $T$, we have
the equivalences:
$T$ is nonexpansive
$\Leftrightarrow$
$(\forall x \in X)$ $\|Tx\|\leq\|x\|$
$\Leftrightarrow$
$(\forall x\in X)$ $\|Tx\|^2\leq\|x\|^2$
$\Leftrightarrow$
$(\forall x\in X)$ $\scal{Tx}{Tx}\leq\scal{x}{x}$
$\Leftrightarrow$
$(\forall x\in X)$ $\scal{x}{T^* Tx}\leq\scal{x}{x}$
$\Leftrightarrow$
$(\forall x\in X)$ $0\leq \scal{x}{(\Id-T^*T)x}$
$\Leftrightarrow$
$\Id-T^* T$ is positive semidefinite. 
\end{proof}

\cref{l:nonexpchar} is needed for the following characterization: 

\begin{proposition}
\label{p:prematrix}
Let $T\colon X\to X$ be linear and nonexpansive, and let $\kappa\in\left]0,1\right]$. 
Then $T$ is $\kappa$-averaged if and only if 
\begin{equation}
(2\kappa-1)\Id - \big(T^*T-(1-\kappa)(T+T^*)\big) \;\;\text{is positive semidefinite.}
\end{equation}
\end{proposition}
\begin{proof}
When $\kappa=1$, this is precisely \cref{l:nonexpchar}.
Using \cite[Proposition~4.35]{BC2017} and \cref{l:nonexpchar}, we have:
\begin{subequations}
\begin{align}
\text{$T$ is $\kappa$-averaged}
&\Leftrightarrow 
\tfrac{1}{\kappa}\big(T-(1-\kappa)\Id\big)
\;\;\text{is nonexpansive}\\
 &\Leftrightarrow 
 \Id - \tfrac{1}{\kappa^2}\big(T-(1-\kappa)\Id\big)^*\big(T-(1-\kappa)\Id\big)
 \succeq 0\\
 &
 \Leftrightarrow 
 \kappa^2\Id - \big(T-(1-\kappa)\Id\big)^*\big(T-(1-\kappa)\Id\big)
 \succeq 0\\
&\Leftrightarrow 
\kappa^2\Id-\big(T^*T-(1-\kappa)(T+T^*)+(1-\kappa)^2\Id \big)\succeq 0\\
&\Leftrightarrow 
(2\kappa-1)\Id-\big(T^*T-(1-\kappa)(T+T^*) \big)\succeq 0,
\end{align}
\end{subequations}
as claimed. 
\end{proof}

We conclude this subsection by providing a procedure 
that allows us to (at least numerically) compute
the modulus of averagedness of a nonexpansive square matrix.

\begin{theorem}
\label{t:onematrix}
Let $M\in\RR^{n\times n}$ be a square matrix that corresponds to
a nonexpansive operator on $\RR^n$.
Set 
$A := (\Id-M){\tran}(\Id-M)$
and 
$B := 2\Id-(M+M{\tran})$.
Suppose that $U:=[u_1|u_2|\ldots|u_d]$, where 
$u_1,\ldots,u_d$ is an orthonormal list of eigenvectors corresponding to all positive eigenvalues $\beta_1,\ldots,\beta_d$ of $B$, ordered (WLOG) decreasingly. 
If $d=0$, then $\kappa(M)=0$. Otherwise,
create the $d\times d$ invertible diagonal matrix $D$ with diagonal elements
$\sqrt{\beta_1},\ldots,\sqrt{\beta_d}$.
Then 
\begin{equation}
\kappa(M) = \text{the largest eigenvalue of}\;\; D^{-1}U{\tran} AUD^{-1}.
\end{equation}
	
\end{theorem}
\begin{proof}
Let $\kappa\in[0,1]$.
It follows from 
\cref{p:prematrix} 
that 
\begin{subequations}
\label{e:230310a}
\begin{align}
&\qquad\text{$M$ is $\kappa$-averaged}
\nonumber
\\
&\siff\ (2\kappa-1)\Id - \big(M\tran M-(1-\kappa)(M+M\tran)\big) 
\succeq 0
\\
&\siff\
2\kappa\Id-\Id-M\tran M+(M+M\tran)-\kappa(M+M\tran)\succeq 0
\\
&\siff\
\kappa (2\Id-(M+M{\tran}))\succeq \Id+M\tran M-(M+M\tran)=(\Id-M){\tran}(\Id-M)
\\
&\siff\
\kappa B\succeq A.
\end{align}
\end{subequations}
Because $M$ is nonexpansive, we observe that \cref{e:230310a} holds with $\kappa=1$. 
Clearly, to find $\kappa(M)$, we must find the smallest nonnegative $\kappa$ making 
\cref{e:230310a} true. 
The result now follows from  \cref{p:superRayleigh}.
\end{proof}

\subsection{A nonlinear result}

\begin{proposition}
Let $g\colon\RR\to\RR$ be a differentiable function. 
Then $g$ is nonexpansive if and only if $|g'|\leq 1$ in which case 
its modulus of averagedness $\kappa(g) = (1-\inf g'(\RR))/2$. 
\end{proposition}
\begin{proof}
The characterization of nonexpansiveness follows easily from the
Mean Value Theorem.
Now assume $g$ is nonexpansive and let $\kappa\in\left]0,1\right]$.
Then $-1\leq\mu\leq \sup g'(\RR)\leq 1$.
Hence we have the equivalences
$g$ is $\kappa$-averaged
$\Leftrightarrow$
$(g-(1-\kappa)\Id)/\kappa$ is nonexpansive
$\Leftrightarrow$
$|(g'-(1-\kappa))/\kappa|\leq 1$
$\Leftrightarrow$
$|g'-(1-\kappa)|\leq \kappa$
$\Leftrightarrow$
$-\kappa\leq g'-(1-\kappa)\leq \kappa$
$\Leftrightarrow$
$1-2\kappa\leq g'\leq 1$
$\Leftrightarrow$
$1-2\kappa\leq \mu$
$\Leftrightarrow$
$(1-\mu)/2 \leq\kappa$. 
\end{proof}

\section{Main result}

\label{sec:main}

We are now ready for our main result. We start with a setting that allows to easily obtain the conclusion: 

\begin{lemma}
\label{l:premain}
Let $U,V$ be two closed linear subspaces of $X$, 
let $\beta\in\zeroun$, and 
set 
\begin{equation}
T := \big((1-\beta)\Id+\beta R_V\big)P_U.
\end{equation}
Suppose that $U=X$ or $U\subseteq V$. 
Then exactly one of the following holds:
\begin{enumerate}
\item 
\label{l:premain1} 
$U=V=X$ and $\kappa(T)=0$.
\item 
\label{l:premain2} 
$U=X$, $V\neq X$, and $\kappa(T)=\beta$. 
\item 
\label{l:premain3} 
$U\neq X$ and $U\subseteq V$, and $\kappa(T)=\thalb$. 
\end{enumerate}
\end{lemma}
\begin{proof}
\cref{l:premain1}:
If $U=V=X$, then $T=\Id$ and thus $\kappa(T)=0$.

\cref{l:premain2}:
Now assume that $U=X$ and $V\neq X$. 
Then $T=(1-\beta)\Id+\beta R_V$. 
It thus follows from \cref{ex:sarah} that $\kappa(T)=\beta$.

\cref{l:premain3}: 
Finally assume that $U\neq X$ and $U\subseteq V$. 
Then $P_VP_U=P_U$ and hence 
\begin{subequations}
\begin{align}
T &= \big((1-\beta)\Id+\beta R_V\big)P_U
= \big((1-\beta)\Id+\beta(2P_V-\Id)\big)P_U\\
&= \big((1-2\beta)\Id+2\beta P_V)\big)P_U
= (1-2\beta)P_U+2\beta P_VP_U\\
&= P_U
= \thalb\Id+\thalb R_U.
\end{align}
\end{subequations}
Again using \cref{ex:sarah}, we see that $\kappa(T)=\thalb$.
\end{proof}

To formulate the main result, we recall that 
the cosine of the \emph{Dixmier} (also known as \emph{minimal}) \emph{angle} of
two closed linear subspaces $U,V$ of $X$ is 
\begin{align}
c_D(U,V) &:= 
\sup\menge{\tscal{u}{v}}{u\in U,\, v\in V,\,\|u\|\leq 1,\,\|v\|\leq 1} 
\end{align}
while the cosine of the \emph{Friedrichs angle} is 
\begin{align}
\label{e:c_F}
c_F(U,V) &:= 
\sup\menge{\tscal{u}{v}}{u\in U\cap(U\cap V)^\perp, \, v\in V\cap(U\cap V)^\perp,\,\|u\|\leq 1,\,\|v\|\leq 1}. 
\end{align}
See \cite{Deutsch} and \cite{Maratea} for (much) more on these notions. 

\begin{theorem}[main result]
\label{t:main}
Let $U,V$ be two closed linear subspaces of $X$, 
let $\beta\in\zeroun$, and set 
\begin{equation}
T := \big((1-\beta)\Id+\beta R_V\big)P_U.
\end{equation}
Suppose that 
$U\subsetneqq X$ and $U\not\subseteq V$. 
Then 
\begin{equation}
\label{e:main}
\kappa(T) = 
\frac{1+2\beta(1-c_F^2)+\sqrt{(1-2\beta)^2+4(1-\beta)\beta c_F^2}}{2(2-\beta c_F^2)}, 
\end{equation}
where $c_F := c_F(U,V)$ denotes the cosine of the Friedrichs angle between $U$ and $V$.
\end{theorem}
\begin{proof}
Recall that $R_V = P_V-P_{V^\perp}$ and so $R_V-\Id = P_V-P_{V^\perp}-(P_V+P_{V^\perp})= -2P_{V^\perp}$. 
Hence 
\begin{subequations}
\begin{align}
T &= \big(\Id +\beta(R_V-\Id)\big)P_U
=\big(\Id-2\beta P_{V^\perp}\big)\big(\Id-P_{U^\perp}\big)\\
&= \Id-P_{U^\perp}-2\beta P_{V^\perp} + 2\beta P_{V^\perp}P_{U^\perp}
\end{align}
\end{subequations}
and thus
\begin{equation}
\Id-T = P_{U^\perp}+2\beta P_{V^\perp} - 2\beta P_{V^\perp}P_{U^\perp}.
\end{equation}
By \cite[Corollary~4.51]{BC2017}, $\Fix T = U \cap V$.
Because $U\neq X$, it follows that 
$\Fix T = U\cap V\subsetneqq X$.
Let $z \in X\smallsetminus \Fix T$  and 
split $z$ into 
\begin{equation}
z=  u+u^\perp,\quad \text{where $u=P_Uz$ and $u^\perp = P_{U^\perp}z$.}
\end{equation}
Then 
\begin{subequations}
\label{e:230304a}
\begin{align}
z-Tz
&= P_{U^\perp}z+2\beta P_{V^\perp}z - 2\beta P_{V^\perp}P_{U^\perp}z
=u^\perp +2\beta P_{V^\perp}u+2\beta P_{V^\perp}u^\perp
-2\beta P_{V^\perp}u^\perp\\
&= u^\perp + 2\beta P_{V^\perp}u.
\end{align}
\end{subequations}
For future reference, we note that  $z\neq Tz$ yields
\begin{equation}
\label{e:230301a}
u^\perp \neq -2\beta P_{V^\perp}u.
\end{equation}
From \cref{e:230304a} we get 
\begin{align}
\|z-Tz\|^2 
&= 
\|u^\perp\|^2 + 4\beta^2\|P_{V^\perp}u\|^2
+4\beta\tscal{u^\perp}{P_{V^\perp}u}
\end{align}
and also  
\begin{subequations}
\begin{align}
\scal{z}{z-Tz}
&=\tscal{u+u^\perp}{u^\perp + 2\beta P_{V^\perp}u}
=\|u^\perp\|^2 + 2\beta\tscal{u+u^\perp}{P_{V^\perp}u} \\
&=\|u^\perp\|^2
+2\beta\|P_{V^\perp}u\|^2 + 2\beta\tscal{u^\perp}{P_{V^\perp}u}.
\end{align}
\end{subequations}
Now write
\begin{equation}
\label{e:thec}
\tscal{u^\perp}{P_{V^\perp}u}=c\|u^\perp\|\|P_{V^\perp}u\|,
\end{equation}
where $c\in [-1,1]$ and where we make $c$ unique 
by setting it equal to $0$ if $u^\perp=0$ or 
$P_{V^\perp}u=0$. 
We have $z\neq Tz$, which  in view of \cref{l:moachar},
makes the following quotient is well defined:
\begin{subequations}
\label{e:theq}
\begin{align}
\frac{1}{2}\frac{\|z-Tz\|^2}{\scal{z}{z-Tz}}
&= \frac{1}{2}\frac{\|u^\perp\|^2 + 4\beta^2\|P_{V^\perp}u\|^2
+4\beta\tscal{u^\perp}{P_{V^\perp}u}}{\|u^\perp\|^2
+2\beta\|P_{V^\perp}u\|^2 + 2\beta\tscal{u^\perp}{P_{V^\perp}u}}\\
&= \frac{\|u^\perp\|^2 + 4\beta^2\|P_{V^\perp}u\|^2
+4\beta c\|u^\perp\|\|P_{V^\perp}u\|}{2\|u^\perp\|^2
+4\beta\|P_{V^\perp}u\|^2 + 4\beta c\|u^\perp\|\|P_{V^\perp}u\|}. 
\end{align}
\end{subequations}
In view of \cref{l:moachar}, our goal is to \emph{maximize} this quotient. 
In view of \cref{e:230301a}, it is impossible that $u^\perp$ and 
$P_{V^\perp}u$ are both simultaneously equal to $0$. 
We now discuss cases.

\emph{Case~1:} $u^\perp=0$.\\ 
Then $P_{V^\perp} u\neq 0$.
Then the quotient in \cref{e:theq} simplifies to
\begin{equation}
\frac{1}{2}\frac{\|z-Tz\|^2}{\scal{z}{z-Tz}} = \beta. 
\end{equation}

\emph{Case~2:} $P_{V^\perp}u=0$.\\
Then $u^\perp\neq 0$. 
Here 
the quotient in \cref{e:theq} simplifies to
\begin{equation}
\frac{1}{2}\frac{\|z-Tz\|^2}{\scal{z}{z-Tz}} = \frac{1}{2}. 
\end{equation}

\emph{Case~3:} $u^\perp\neq 0$ and $P_{V^\perp}u\neq 0$.\\ 
In this case, $c$ from \cref{e:thec} is is equal to the cosine between the 
vectors $u^\perp$ and $P_{V^\perp}u\neq 0$:
\begin{equation}
c = \frac{\tscal{u^\perp}{P_{V^\perp}u}}{\|u^\perp\|\|P_{V^\perp}u\|}. 
\end{equation}
We rewrite the quotient from \cref{e:theq} as
\begin{subequations}
\label{e:zquotient}
\begin{align}
\frac{1}{2}\frac{\|z-Tz\|^2}{\scal{z}{z-Tz}}
&= \frac{\|u^\perp\|^2 + 4\beta^2\|P_{V^\perp}u\|^2
+4\beta c\|u^\perp\|\|P_{V^\perp}u\|}{2\|u^\perp\|^2
+4\beta\|P_{V^\perp}u\|^2 + 4\beta c\|u^\perp\|\|P_{V^\perp}u\|}\\
&= \frac{(\|u^\perp\|/\|P_{V^\perp}u\|)^2 + 4\beta^2
+4\beta c(\|u^\perp\|/\|P_{V^\perp}u\|)}{2(\|u^\perp\|/\|P_{V^\perp}u\|)^2
+4\beta + 4\beta c(\|u^\perp\|/\|P_{V^\perp}u\|)}\\
&= \frac{t^2 + 4\beta ct + 4\beta^2}{2t^2 + 4\beta c t + 4\beta},
\quad \text{where $t= \|u^\perp\|/\|P_{V^\perp}u\|>0$.}
\end{align}
\end{subequations}
Because we can vary $u$ and $u^\perp$ independently, we can arrange for $t$ to be
any positive number.
We now discuss this quotient as a function of $t$:
\begin{equation}
q(t) := 
\frac{t^2 + 4\beta ct + 4\beta^2}{2t^2 + 4\beta c t + 4\beta}
\end{equation}
where $t>0$. Clearly, $q$ is continuous (even analytic). 
Note that as $t\to 0^+$, we recover \emph{Case~1} while
$t\to\pinf$ leads to \emph{Case~2}. 
For convenience, we set $q(0) := \lim_{t\to 0^+} q(t)=\beta$ and 
$q(+\infty) = \tfrac{1}{2}$. 
So \emph{Case~3} is the only case we need to focus on. 
We shall also sometimes view $q(t)$ as a function in $c$ as well --- 
in that case, we emphasize this by writing $q(t,c)$.
Next, 
\begin{align}
\label{e:230302f}
q'(t)
&=
\frac{\beta \big( -c t^2 + (2-4\beta)t + 4\beta(1-\beta)c\big)}{(t^2 + 2\beta c t + 2\beta )^2}.
\end{align}
Thus 
\begin{equation}
\label{e:230302d}
q'(t)=0
\;\;\Leftrightarrow\;\;
-ct^2+(2-4\beta)t + 4\beta(1-\beta)c=0.
\end{equation}
To make progress, we split our work into additional cases. 

\emph{Case~3.1:} $c=0$.\\
If $\beta=\thalb$, then $q'\equiv 0$ and $q$ is the constant function $\thalb$:
\begin{equation}
\label{e:230302a}
\beta=\thalb \;\;\Rightarrow\;\; \sup_{t>0} q(t,0)=\thalb. 
\end{equation}
Now assume that $\beta\neq\thalb$. Then $q'(t)=0$ $\Leftrightarrow$ $t=0$.
Thus, $q$ has no strictly positive critical point. It follows that 
$q(\cdot,0)$ is either increasing or decreasing on $\RPP$. 
Hence 
\begin{equation}
\label{e:230302b}
\beta\neq \thalb \;\;\Rightarrow\;\; \sup_{t>0} q(t,0)=\max\{q(0,0),q(\pinf,0)\}
=\max\{\beta,\thalb\}. 
\end{equation}
Because the formula from \cref{e:230302b} still holds in \cref{e:230302a},
we combine and summarize \emph{Case~3.1} into
\begin{equation}
\label{e:230302bb}
\sup_{t>0} q(t,0)=\max\{q(0,0),q(\pinf,0)\}
=\max\{\beta,\thalb\}.
\end{equation}

\emph{Case~3.2:} $c\neq 0$.\\
In view of \cref{e:230302d}, $q'(t)$ has two real roots, namely
\begin{equation}
\label{e:230302e}
\frac{1-2\beta \pm \sqrt{(1-2\beta)^2+4(1-\beta)\beta c^2}}{c}. 
\end{equation}
We care about the positive root. Because this depends on the sign of $c$,
we bifurcate the analysis further. 

\emph{Case~3.2.1:} $c> 0$.\\
It follows from \cref{e:230302e} that the unique positive root of $q'$ is 
\begin{equation}
t_1 := \frac{1-2\beta + \sqrt{(1-2\beta)^2+4(1-\beta)\beta c^2}}{c}. 
\end{equation}
On the other hand, \cref{e:230302f} shows that for all $t$ sufficiently large,
we have $q'(t)<0$. This shows that left of $t_1$, the function $q'$ is positive
while to the right of $t_1$, the function $q'$ is negative. 
Hence $t_1$ is a \emph{local maximizer} of $q$. Because $q$ has no further positive critical point,
it is clear that $t_1$ is the global maximizer of $q|_{\RP}$:
$q(t_1)=\max_{t\geq 0}q(t)$. Computing $q(t_1)$ and simplifying leads to 
\begin{equation}
\max_{t > 0}q(t) = \frac{1+2\beta(1-c^2)+\sqrt{(1-2\beta)^2+4(1-\beta)\beta c^2}}{2(2-\beta c^2)} > \max\{q(0),q(+\infty)\}= \max\{\beta,\thalb\}.
\end{equation}

\emph{Case~3.2.2:} $c< 0$.\\
It follows this time from \cref{e:230302e} that the unique positive root of $q'$ is 
\begin{equation}
t_2 := \frac{1-2\beta - \sqrt{(1-2\beta)^2+4(1-\beta)\beta c^2}}{c}. 
\end{equation}
On the other hand, \cref{e:230302f} shows that for all $t$ sufficiently large,
we have $q'(t)>0$. This shows that left of $t_2$, the function $q'$ is negative 
while to the right of $t_2$ the function $q'$ is positive.
Hence $t_2$ is a \emph{local minimizer} of $q$. 
Because $q$ has no further positive critical point,
it is clear that $t_2$ is the minimizer of $q|_{\RP}$ and that 
$q(t_2)<\min\{q(0),q(\infty)\}\leq \max\{q(0),q(\infty)\}= \max\{\beta,\thalb\}$.
So the critical point $t_2$ does not help in the maximization of $q(t)$. 
Finally note that if $c<0$, then we can reach \emph{Case~3.2.1} by simply replacing 
either $u^\perp$ or $u$ (but not both) by its negative. 

We thus sum up \emph{Case~3.2} as follows:
\begin{subequations}
\begin{align}
\max\{ \sup_{t>0} q(t,c),\sup_{t>0} q(t,-c)\} 
&= \frac{1+2\beta(1-c^2)+\sqrt{(1-2\beta)^2+4(1-\beta)\beta c^2}}{2(2-\beta c^2)}
\label{e:230303a}\\
&> \max\{q(0),q(\infty)\}= \max\{\beta,\thalb\}. 
\end{align}
\end{subequations}
Even though $c\neq 0$ in this current \emph{Case~3.2} under consideration, 
let us replace 
$c$ in the right side of \cref{e:230303a} by $c=0$.
We then obtain 
\begin{align}
\frac{1+2\beta +\sqrt{(1-2\beta)^2}}{4}
&= \frac{1+2\beta+|1-2\beta|}{4}
= \frac{\thalb+\beta+|\thalb-\beta|}{2}
=\max\{\thalb,\beta\},
\end{align}
which we saw already in \cref{e:230302bb}! 
This allows us to \emph{combine} \emph{Case~3.1} and \emph{Case~3.2} and record
a single formula valid for \emph{Case~3}, namely:
\begin{equation}
\sup_{t>0} q(t,\pm c) = 
\frac{1+2\beta(1-c^2)+\sqrt{(1-2\beta)^2+4(1-\beta)\beta c^2}}{2(2-\beta c^2)}. 
\end{equation}
Having resolved the maximization over $t$, it remains to maximize over $c\in [0,1]$. 
To this end, we study the continuous function
\begin{equation}
Q\colon [0,1]\to\RPP\colon x \mapsto 
\frac{1+2\beta(1-x^2)+\sqrt{(1-2\beta)^2+4(1-\beta)\beta x^2}}{2(2-\beta x^2)}. 
\end{equation}
It turns out that on $\zeroun$, we have 
\begin{equation}
\frac{dQ}{dx} = 
\frac{\beta x}{r(2-\beta x^2)^2}\big(n_1+n_2+n_3\big), 
\end{equation}
where 
\begin{subequations}
\begin{align}
r &:= \sqrt{(1-2\beta)^2+4(1-\beta)\beta x^2} > 0,\\ 
n_1 &:= 2\beta(1-\beta) x^2 > 0,\\
n_2 &:= 4(1-\beta)^2+1 > 0,\\
n_3 &:= (2\beta-3)r < 0.  
\end{align}
\end{subequations}
We now make the following

\emph{Claim:} $n_1+n_2+n_3>0$.\\
The \emph{Claim} is equivalent to $n_1+n_2>-n_3 = |n_3|$ and hence to 
$(n_1+n_2)^2 > n_3^2$. Luckily, we have
\begin{equation}
(n_1+n_2)^2-n_3^2 = 4(1-\beta)^2(2-\beta x^2)^2>0. 
\end{equation}
This verifies the \emph{Claim} which in turn yields that 
\begin{equation}
\label{e:Qincreases}
\text{$Q$ is strictly increasing.}
\end{equation}
This implies that 
\begin{equation}
\kappa(T) = Q(c_D),
\quad\text{where~}
c_D := c_D\big(U^\perp,\overline{P_{V^\perp}(U)}\big). 
\end{equation}
Luckily, this Dixmier angle can be simplified further.
Indeed, 
using \cite[Corollary~15.31(i)]{BC2017}, we have 
\begin{align}
\overline{P_{V^\perp}(U)}
&=\big(P_{V^\perp}^{-1}(U^\perp)\big)^\perp
=
\big((U^\perp\cap V^\perp)+V \big)^\perp
=
\big(U^\perp\cap V^\perp)^\perp \cap V^\perp.
\end{align}
Hence
\begin{equation}
c_D = c_D\big(U^\perp, V^\perp \cap (U^\perp\cap V^\perp)^\perp  \big). 
\end{equation}
On the other hand, \cite[Lemma~10.1 and Theorem~16]{Maratea} yield that 
\begin{equation}
c_D(U^\perp, V^\perp \cap (U^\perp\cap V^\perp)^\perp  )= c_F(U^\perp,V^\perp)
= c_F(U,V) =: c_F. 
\end{equation}
Altogether, we conclude that 
\begin{equation}
\kappa(T) = Q(c_F),
\end{equation}
and we are done.
\end{proof}

A very pleasing formula arises when $\beta=\thalb$ in the previous results.

\begin{corollary}
\label{c:bonito}
Let $U,V$ be closed linear subspaces of $X$. 
Then either $U=V=X$ and $\kappa(P_VP_U)=0$, or 
\begin{equation}
\label{e:bonito}
\kappa(P_VP_U) = \frac{1+c_F} {2+c_F}, \quad
\text{where~}
c_F=c_F(U,V). 
\end{equation}
\end{corollary}
\begin{proof}
This follows from \cref{l:premain} and \cref{t:main}, after noting
(see \cite[Lemma~9.5]{Deutsch}) that 
$c_F = 0$ when $U\subseteq V$ or $V\subseteq U$.  
\end{proof}

\begin{remark}[Ogura--Yamada bound is optimal]
In the setting of \cref{t:main},
the operator $T$ is the composition of two averaged maps:
$P_U$ is $\thalb$-averaged and 
$(1-\beta)\Id+\beta R_V$ is $\beta$-averaged (see \cref{ex:sarah}). 
Hence the Ogura--Yamada bound (see \cref{e:OY}) in this case turns into 
\begin{equation}
\label{e:OY1}
\frac{\thalb + \beta -2\thalb\beta}{1-\thalb \beta}
= \frac{1}{2-\beta}
\end{equation}
which coincides precisely with $\kappa(T)$ from \cref{e:main} assuming we $c_F=1$.
On the one hand, if $X$ is finite-dimensional, then 
$c_F(U,V)<1$ (see, e.g., \cite[Theorem~13]{Maratea}) and hence 
$\kappa(T)$ from \cref{e:main} is sharper than the Ogura--Yamada bound
from \cref{e:OY1}. 
On the other hand, if $X$ is infinite-dimensional, then there exist
incarnations of two disjoint subspaces such that $c_F=1$ 
(see, e.g., \cite[Section~6]{86}). 
This shows
that the Ogura--Yamada bound is optimal!
\end{remark}

Let us now revisit an example used to show that 
the composition of two firmly nonexpansive mappings may fail to be firmly nonexpansive:

\begin{example}
\label{ex:wow}
Suppose that $X=\RR^2$ and consider the matrix
\begin{equation}
M := \thalb \begin{bmatrix}
1 & 0\\1 & 0
\end{bmatrix}, 
\end{equation}
which represents $P_VP_U$, where 
$U=\RR[1,0]{\tran}$ and $V=\RR[1,1]{\tran}$. 
It is pointed out in \cite[Example~4.45]{BC2017} that $M$ is 
$\tfrac{2}{3}$-averaged but not $\thalb$-averaged. 
It is clear that the Friedrichs (and Dixmier) angle between $U$ and $V$
is $45^\circ = \pi/4$. 
Hence $c_F := c_F(U,V) = \cos(\pi/4) = 1/\sqrt{2}$.
Applying \cref{c:bonito}, we thus learn that 
the true modulus of $M$ is 
\begin{equation}
\kappa(M) = \frac{1+1/\sqrt{2}}{2+1/\sqrt{2}} = \frac{1+\sqrt{2}}{1+2\sqrt{2}}
= \frac{3+\sqrt{2}}{7}\approx 0.6306 < 0.\overline{6}.
\end{equation}
A derivation of this modulus via \cref{t:onematrix} is presented in \cref{app2}. 
Because we now know the answer, we also provide an elementary \emph{a posteriori} verification 
$\kappa(M)$ in \cref{app3}. 
\end{example}

\begin{remark}
Because of \cref{p:mod*}, we could also write down a result for $\kappa(T)$ 
when $T = P_U((1-\beta)\Id+\beta R_V)$. 
\end{remark}

\begin{remark}[future work]
It would be desirable to obtain $\kappa(T)$ when 
$T = ((1-\beta)\Id+\beta R_V)((1-\alpha)\Id+\alpha R_U)$, which would 
subsume the present work when $\alpha=1/2$. 
Unfortunately, the current proof breaks down in this setting and does not lead
to an exact result. 
We do conjecture that the Ogura--Yamada bound from \cref{e:OY} 
is sharp in this case as well. 
\end{remark}

\section*{Acknowledgments}
\small
The research of the authors was partially supported by Discovery Grants
of the Natural Sciences and Engineering Research Council of
Canada.

\appendix 

\section{A result utilized in the proof of \cref{t:onematrix}}

\label{app1}

We suspect the following result is well known; however, we were not able to find
a clear reference so we provide a self-contained proof. 

\begin{proposition}
\label{p:superRayleigh}
Let $A,B\in\RR^{n\times n}$ be symmetric and positive semidefinite.
Suppose that there exists  $\beta\geq 0$ such that $A\preceq\beta B$.
Denote the nonnegative eigenvectors of $B$ by 
$\beta_1\geq \beta_2\geq \cdots \geq \beta_d>\beta_{d+1}=\cdots=\beta_n=0$.
If $d=0$, then set $\gamma=0$. 
If $d\geq 1$, then proceed as follows:
Obtain orthonormal eigenvectors $u_1,\ldots,u_d$ of $U$ with corresponding
eigenvalues $\beta_1,\ldots,\beta_d$.
Build the matrix $U:=[u_1 |u_2 |\ldots | u_d]\in\RR^{n\times d}$, 
the diagonal matrix $D\in \RR^{d\times d}$ with diagonal entries 
$\sqrt{\beta_1},\ldots,\sqrt{\beta_d}$, and finally the matrix 
$C := D^{-1}U{\tran} AUD^{-1}$. 
Denote its largest eigenvalue by $\gamma$.
No matter which case we are in ($d=0$ or $d\geq 1$), the $\gamma$ constructed
is the smallest nonnegative real number such that 
$A\preceq \gamma B$.
\end{proposition}
\begin{proof}
The assumptions $0\preceq A\preceq \beta B$ imply that 
\begin{equation}
\label{e:assumpAB}
(\forall x\in\RR^n)\quad 
0 \leq \scal{x}{Ax}\leq\beta\scal{x}{Bx}. 
\end{equation}
We start by proving the following claim:
\begin{equation}
\label{e:ker:inc}
\ker B = \menge{x\in \RR^n}{\scal{x}{Bx}=0}
\subseteq \menge{x\in \RR^n}{\scal{x}{Ax}=0} = \ker A.
\end{equation}
Note that if $x\in\ker B$, then $\scal{x}{Bx}=\scal{x}{0}=0$ ---
the converse also holds as we show next:
Consider spectral decomposition  
$B = \sum_{i=1}^n \beta_iP_i$, 
where  
$P_i := u_iu_i{\tran}$ and $u_1,\ldots,u_n$ is a list of orthonormal eigenvectors 
corresponding to the eigenvalues $\beta_1,\ldots,\beta_n$ of $B$. 
We note that if $d=0$, then $B=0$ and the conclusion is clear.
So assume that $d\geq 1$ and let
$x\in\RR^n$ satisfy $\scal{x}{Bx}=0$.
Then $0=\tscal{x}{{\sum_{i\leq d}\beta_iP_ix}}=\sum_{i\leq d}\beta_i\scal{x}{P_ix}
=\sum_{i\leq d}\beta_i\|P_ix\|^2$. Hence 
if $i\leq d$, then $P_ix=0$. 
It follows that $x=\sum_{i}P_ix = \sum_{i>d}P_ix\in\ker B$.
We have shown
\begin{equation}
\ker B = \menge{x\in \RR^n}{\scal{x}{Bx}=0}.
\end{equation}
We have analogously that $\ker A =\menge{x\in \RR^n}{\scal{x}{Bx}=0}$. 
Moreover, 
\cref{e:assumpAB} implies that 
$\scal{x}{Bx}=0$ $\Rightarrow$ $\scal{x}{Ax}=0$. 
This proves \cref{e:ker:inc}.

Next, we observe that because $A,B$ are symmetric matrices, we also have
$\ran A \perp \ker A$ and $\ran B \perp \ker B$. 
Let's decompose $\RR^n$ into $\RR^n = K\oplus R$, where $K=\ker B$ and 
$R = \ran B$. Write $x=k+ r$, where $k=P_Kx$, $r=P_Rx$, and $k\perp r$.
Then 
$\scal{x}{Bx}=\scal{x}{Bk}+\scal{x}{Br}=\scal{x}{0}+\scal{x}{Br}=
\scal{k+r}{Br}=\scal{k}{Br}+\scal{r}{Br}=\scal{r}{Br}$
where the last equality follows because $k\in K$, $Br\in\ R$, and $K\perp R$.  
Hence 
\begin{equation}
\scal{x}{Bx}=\scal{r}{Br}.
\end{equation}
Similarly,
$\scal{x}{Ax}=\scal{x}{Ak}+\scal{x}{Ar}=\scal{x}{Ar}$
because $k\in K = \ker B \subseteq \ker A$ by \cref{e:ker:inc} and so $Ak=0$.
Furthermore,
$\scal{x}{Ar}=\scal{k}{Ar}+\scal{r}{Ar}=\scal{r}{Ar}$ 
because $k\in\ker A$, $Ar\in\ran A$, and $\ker A\perp \ran A$. 
Hence we also
\begin{equation}
\scal{x}{Ax}=\scal{r}{Ar}.
\end{equation}
Thus the problem of finding the optimal $\gamma$ is reduced to
\begin{equation}
(\forall r\in R\smallsetminus\{0\})\quad 
0 \leq \scal{r}{Ar}\leq\gamma\scal{r}{Br}. 
\end{equation}
But $r \in R\smallsetminus\{0\}$ means also that $r\notin \ker B$ and so $\scal{r}{Br}>0$.
Dividing now reduces this to 
\begin{equation}
(\forall r\in R\smallsetminus\{0\})\quad 
\frac{\scal{r}{Ar}}{\scal{r}{Br}}\leq\gamma. 
\end{equation}
So finding $\gamma$ means maximizing the quotient and we also may assume that 
$\|r\|=1$.

Next, we note that
$R=\ran B = \spn\{u_1,\ldots,u_d\}=\ran U$. 
Any $r\in R$ can be written as $r=Uy$, where $\|y\|=\|r\|$ and $y\in\RR^d$.
Moreover,
$\scal{r}{Br}=r{\tran} Br=y{\tran} U{\tran} BUy
= \tscal{y}{U{\tran} BU y}$ 
and $\scal{r}{Ar}=r{\tran} Ar=y{\tran} U{\tran} AUy
= \tscal{y}{U{\tran} AU y}$.
Hence we are led to the quotient
\begin{equation}
(\forall y\in \RR^d\smallsetminus\{0\})\quad
\frac{\tscal{y}{U{\tran} AU y}}{\scal{y}{U{\tran} BU y}}\leq\gamma.
\end{equation}
The matrices involved are of size $d\times d$, the $U{\tran} BU$ is invertible, and
by construction of $D$, we have 
$U{\tran} BU = D^2$. 
Note that 
\begin{equation}
C = D^{-1}U{\tran} AUD^{-1} 
\end{equation}
is a symmetric matrix. 
Finally, we change variables to $z=D y$.
Then 
\begin{subequations}
\begin{align}
\frac{\tscal{y}{U{\tran} AU y}}{\scal{y}{U{\tran} BU y}}
&=
\frac{\tscal{y}{DD^{-1}U{\tran} AU y}}{\scal{y}{D^2 y}}\\
&=
\frac{\tscal{D y}{D^{-1}U{\tran} AU D^{-1}D y}}{\scal{D y}{D y}}\\
&= \frac{\tscal{z}{D^{-1}U{\tran} AU D^{-1}z}}{\scal{z}{z}}\\
&=
\frac{\scal{z}{Cz}}{\|z\|^2}
\end{align}
\end{subequations}
and the maximum of this quotient is the largest eigenvalue of $C$ by the 
classical Rayleigh theorem \cite[Theorem~4.21]{Bisgard}. 
\end{proof}

\section{Alternative derivation of $\kappa(M)$ from \cref{ex:wow}}

\label{app2}
We compute 
\begin{equation}
A = (\Id-M{\tran})(\Id-M)= 
\frac{1}{2}\begin{bmatrix}
1 & -1\\ -1& 2
\end{bmatrix} 
\quad\text{and}\quad
B = 2\Id-(M+M{\tran}) = \frac{1}{2}\begin{bmatrix}
2 & -1\\ -1& 4
\end{bmatrix}. 
\end{equation}
The eigenvalues of $B$ are 
$\beta_1 := (3+\sqrt{2})/2$ 
and 
$\beta_2 := (3-\sqrt{2})/2$ 
which are both positive. 
The corresponding non-normalized eigenvectors are 
$[1-\sqrt{2},1]{\tran}$ and $[1+\sqrt{2},1]{\tran}$, and their normalized counterparts are 
\begin{align}
u_1 = \frac{[1-\sqrt{2},1]{\tran}}{\|[1-\sqrt{2},1]{\tran}\|}
= 
\frac{[1-\sqrt{2},1]{\tran}}{\sqrt{4-2\sqrt{2}}}
= 
\frac{[-\sqrt{2-\sqrt{2}},\sqrt{2+\sqrt{2}}]{\tran}}{2}
\end{align}
and 
\begin{align}
u_2 = \frac{[1+\sqrt{2},1]{\tran}}{\|[1+\sqrt{2},1]{\tran}\|}
= 
\frac{[1+\sqrt{2},1]{\tran}}{\sqrt{4+2\sqrt{2}}} 
= 
\frac{[\sqrt{2+\sqrt{2}},\sqrt{2-\sqrt{2}}]{\tran}}{2}. 
\end{align}
The matrix $U$ is therefore given by
\begin{align}
U = 
\frac{1}{2}
\begin{bmatrix}
-\sqrt{2-\sqrt{2}} & {\sqrt{2+\sqrt{2}}}\\
{\sqrt{2+\sqrt{2}}} & \sqrt{2-\sqrt{2}}
\end{bmatrix} = U{\tran}
\end{align}
while 
\begin{align}
D = \begin{bmatrix}
\sqrt{\beta_1} & 0 \\
0 &\sqrt{\beta_2}
\end{bmatrix}
= 
\frac{1}{\sqrt{2}}
\begin{bmatrix}
\sqrt{3+\sqrt{2}} & 0\\
0 & \sqrt{3-\sqrt{2}}
\end{bmatrix} 
\;\;\text{and}\;\;
D^{-1} = {\sqrt{2}}\begin{bmatrix}
\frac{1}{\sqrt{3+\sqrt{2}}} & 0 \\ 0 & \frac{1}{\sqrt{3-\sqrt{2}}}
\end{bmatrix}. 
\end{align}
One computes that 
\begin{align}
D^{-1}U{\tran}A UD^{-1} = 
\frac{1}{28}
\begin{bmatrix}
3(4+\sqrt{2}) & -\sqrt{14}\\ -\sqrt{14} & 3(4-\sqrt{2})
\end{bmatrix}. 
\end{align}
The eigenvalues of this last matrix are $(3\pm\sqrt{2})/7$.
It finally follows from \cref{t:onematrix} that $\kappa(M) = (3+\sqrt{2})/7$ as expected.

\section{\emph{A posterior} verification of $\kappa(M)$ from \cref{ex:wow}}

\label{app3}

Using the linearity of $M$, we get that $M$ is $\kappa$-averaged if and only if
\begin{equation}
\label{e:230309a}
\|Mz\|^2 + (1-2\kappa)\|z\|^2\leq 2(1-\kappa)\scal{z}{Mz}
\end{equation}
for every $z\in \RR^2$. 
Now write $z = [x,y]{\tran}$. 
Then $Mz = \thalb[x,x]{\tran}$. 
Hence
$\|Mz\|^2 = \thalb x^2$, 
$\|z\|^2=x^2+y^2$ and 
$\scal{z}{Mz}=\thalb(x^2+xy)$.  
Thus \cref{e:230309a} turns into 
\begin{equation}
\thalb x^2 + (1-2\kappa)(x^2+y^2)\leq 2(1-\kappa)\thalb(x^2+xy); 
\end{equation}
equivalently, 
\begin{equation}
\label{e:230309b}
\big(\kappa-\thalb\big)x^2 + \big(2\kappa-1\big)y^2 + \big(1-\kappa\big)xy \geq 0.
\end{equation}
Now set
\begin{equation}
\kappa := \frac{3+\sqrt{2}}{7} + \frac{\delta}{14} 
\end{equation}
and multiply \cref{e:230309b} by 14 to get: 
\begin{equation}
\label{e:230309c}
\big(2(3+\sqrt{2})+\delta-7\big)x^2 + 
\big(4(3+\sqrt{2})+2\delta-14\big)y^2 + \big(14-2(3+\sqrt{2})-\delta\big)xy \geq 0.
\end{equation}
In turn, \cref{e:230309c} simplifies to 
\begin{equation}
\label{e:230309d}
\big(2\sqrt{2}-1+\delta\big)x^2 + 
\big(4\sqrt{2}-2+2\delta\big)y^2 + \big(8-2\sqrt{2}-\delta\big)xy \geq 0,
\end{equation}
which we now rewrite as 
\begin{equation}
\label{e:230309e}
\big(2\sqrt{2}-1\big)\big(x+\sqrt{2}y \big)^2 + \delta\big(x^2+2y^2-xy\big)\geq 0.
\end{equation}
This last inequality \cref{e:230309e} is clearly satisfied if $\delta=0$.
Hence $\kappa(M)\leq (3+\sqrt{2})/7$. 
However, if $x = -\sqrt{2}y \neq 0$, then the
left side of \cref{e:230309e}
turns into 
$0+\delta(2y^2 + 2y^2+\sqrt{2}y^2) = (4+\sqrt{2})\delta y^2$ and this is
negative if $\delta<0$. 
Thus $\kappa(M)\geq (3+\sqrt{2})/7$.
Altogether, 
$\kappa(M)= (3+\sqrt{2})/7$.

\end{document}